\documentclass[11pt,a4paper,reqno]{amsart}
\usepackage[T1]{fontenc}
\usepackage{csquotes}
\usepackage{amsmath, amsthm, amsfonts,geometry, bm, amssymb,color}
\usepackage{a4wide,amsfonts,amsmath,latexsym,amssymb,euscript,graphicx,units,mathrsfs,bbold}  \usepackage{multicol} 
\usepackage{comment}
\usepackage{bbm}
\usepackage{theoremref}
 
 \newcommand{\R}{{\mathbb{R}}}

\renewcommand{\d}{\mathrm{d}}
\newcommand{\<}{\langle}
\renewcommand{\>}{\rangle}

   \def\xitilde{\Tilde{\xi}}

    \def\ztilde{\Tilde{z}}

 \def\ytilde{\Tilde{y}}
  \def\talpha{t_{\alpha}}

 \newcommand{\ex}[1]{\E\left[#1\right]}
  \newcommand{\pr}[1]{\P\left(#1\right)}
 \usepackage[shortlabels]{enumitem}
   \usepackage{hyperref}
\hypersetup{
    colorlinks=true,
    linkcolor=blue,
    filecolor=magenta,      
    urlcolor=cyan,
       pdftitle={Overleaf Example},
    pdfpagemode=FullScreen,
    }
\usepackage{kantlipsum}

\def\d{\mathrm d}
\def\TV{\mathrm TV}

\def\HH{\mathfrak{H}}

\def\R{\mathbb{R}}
\def\D{\mathbb{D}}
\def\E{\mathrm{E}}
\def\P{\bm{\mathrm{P}}}

\def\Var{{\mathrm{Var}}}
\def\dom{{\mathrm{Dom}}}

  \newcommand{\e}{\varepsilon}

\def\y{\mathcal{Y}}
  
 \def\H{\mathcal{H}}
  
  \newtheorem{theorem}{Theorem}[section]

\newtheorem{lemma}[theorem]{Lemma}
 \newtheorem{proposition}[theorem]{Proposition}

\numberwithin{equation}{section}
\begin{document}
   
   \title[Convergence of densities for the parabolic Anderson model]{Convergence of densities of spatial averages of the parabolic Anderson model driven by colored noise}

 \author[S. Kuzgun]{Sefika Kuzgun}
\address{University of Kansas, Department of Mathematics, USA}
\email{sefika.kuzgun@ku.edu}

\author[D. Nualart]{David Nualart} \thanks{%
	The work by D. Nualart has been supported   by the  NSF grants DMS-2054735}
\address{University of Kansas, Department of Mathematics, USA}
\email{nualart@ku.edu}

\begin{abstract} 
In this paper, we present a rate of convergence in the uniform norm for the  densities of
 spatial averages of the solution to the $d$-dimensional parabolic
 Anderson model driven by a Gaussian multiplicative noise, which is white in time and has a spatial covariance given by the 
 Riesz kernel.  The proof is based on the combination of
  Malliavin calculus techniques and the Stein’s method for normal approximations.


\end{abstract}

\maketitle


\section{Introduction}
Consider the parabolic Anderson model
 \begin{equation} \label{PAM}
  \frac{\partial u}{\partial t}= \frac{1}{2}  \Delta u+ u \Dot{W},  \qquad x\in \R^d, \,\,   t>0, 
\end{equation} 
with initial condition $u(0,x)=1$, where $\Dot{W}$  is a noise which is white in time and colored in space. This is to say, informally, 
that $\Dot{W}=\{ \Dot{W}(t,x): (t,x) \in \R_+\times \R^d\}$ is a centered Gaussian random field with covariance
\begin{align*}
    \ex{\Dot{W}(t,x)\Dot{W}(s,y)}=\delta_0(t-s)|x-y|^{-\beta},
\end{align*}
for $0<\beta<\min(2,d)$ where $\delta_0$ is the Dirac delta measure at zero.
The existence and uniqueness of a mild solution $u(t,x)$ to  equation \eqref{PAM} has been  proved by Dalang in \cite{Da}, assuming that $u_0$ is bounded. The one-dimensional equation driven by a space-time white noise was studied in the pioneering work by Walsh \cite{Wa}.
Here and  along the paper we will make use of the notation $$\bm{p}_t(x):=\frac{1}{(2\pi t)^{d/2}}e^{-|x|^2/2t}$$ for $t>0$ and $x\in \R^d$. Fix $R>0$ and let $\bm{Q}_R=[-R,R]^d$. Consider the corresponding centered and normalized spatial averages  defined by
\begin{align}
  \label{FR}  F_{R,t}:=\frac{1}{\sigma_{R,t}}\left(\int_{\bm{Q}_R} u(t,x)dx-(2R)^d\right),  \text{ where  } \sigma^2_{R,t}:=\Var\left(\int_{\bm{Q}_R} u(t,x)dx\right). \end{align}

In this paper, we will investigate the quantitative rates of convergence  corresponding to the normal approximations in density of the above spatial averages. The main result of this paper is as follows:

\begin{theorem}
\label{thm:rate1}  Let  $u=\left\{u(t,x): (t,x)\in  \R_+\times \R^d\right\}$ be the mild solution to the stochastic heat equation \eqref{PAM} with the initial condition $u_0=1$. 
Fix $t>0$ and let $F_{R,t}$ be defined as in \eqref{FR}. Then, for all $R\ge1$, $F_{R,t}$ has a density $f_{F_{R,t}}$ and 
\begin{align*}
  \sup_{z\in \R} |f_{F_{R,t}}(z)-\phi(z)| \leq C_t R^{-\beta/2}, 
    \end{align*} 
    where $\phi$ is the density of a standard normal distribution on $\R$.
\end{theorem}

 \section{Preliminaries}

In this section, we will recall some preliminary notions and results. 
Let $C_c^{\infty}(\R_+\times \R^d)$ be the space of infinitely differentiable functions with compact support on $\R_+\times \R^d$. Let $\HH$ be the Hilbert space defined as the completion of $C_c^{\infty}(\R_+\times\R^d)$ with respect to the inner product
\begin{align*}
    \<\varphi,\psi\>_{\HH}=\int_0^{\infty}\int_{\R^{2d}} \varphi(s,x)\psi(s,y)|x-y|^{-\beta}dxdyds.
\end{align*}
Suppose  $W$ a Gaussian noise encoded by a centered Gaussian family of random variables $\{W(\varphi); \varphi \in C_c^{\infty}(\R_+\times \mathbb{R}^d)\}$,  defined in a  complete probability space   $(\Omega, \mathcal{F}, \P)$, with the covariance structure \begin{align} \label{CovarianceSturucture}
   \ex{W(\varphi)W(\psi)}=\<\varphi,\psi\>_{\HH}.
\end{align}
Let $\left(\mathcal{F}_t\right)_{t\in \R_+}$ be the filtration such that $
  \mathcal{F}_t=\mathcal{F}_t^0 \vee \mathcal{N}$
 where $ \mathcal{F}_t^0$ is the $\sigma$-field generated by the random variables $W(\varphi)$ where $\varphi$ has support in $[0,t]\times \mathbb{R}^d$, and $\mathcal{N}$ is the $\sigma$-field generated by the $\P$-null sets.  We say that a random field $X=\{X(t,x): (t,x) \in \mathbb{R}_+ \times \mathbb{R}^d\}$
 is adapted if  for each $(t,x)\in \mathbb{R}_+\times \R^d$ the random variable $X(t,x)$ is $\mathcal{F}_t$-measurable.
Then, for any adapted, jointly measurable random field $X$ such that \begin{align}\label{integrable}
    \ex{\|X\|_{\HH}^2}=\int_0^{\infty}\int_{\R^{2d}} \ex{X(s,x)X(s,y)}|x-y|^{-\beta}dxdyds<\infty,
\end{align}
the stochastic integral 
\begin{align*}
    \int_{\R_+\times\R^d} X(\tau,\xi)W(d\tau,d\xi)
\end{align*}
is well-defined in the sense of Walsh and the It\^o-Walsh isometry 
\begin{align*}
   \ex{\left| \int_{\R_+\times\R^d} X(s,y)W(ds,dy)\right|^2}=\int_0^{\infty}\int_{\R^{2d}} \ex{X(s,x)X(s,y)}|x-y|^{-\beta}dxdyds
\end{align*}
holds. Moreover, the following Burkholder-Davis-Gundy equality is satisfied.

\begin{proposition}\label{prop:BDG} For all $p\geq 2$ there exists a constant $C_p>0$ such that for all adapted, jointly measurable random field $X$ satisfying \eqref{integrable} and for all $t\in [0,\infty)$, we have \begin{align}
    \ex{\left| \int _{[0,t]\times \R^{d}} X(s,y)W(ds,dy)\right|^{p}} \leq C_p \ex{\left(\int_0^t \int _{\R^{2d}} X(s,y)X(s,x)|x-y|^{-\beta}dxdyds\right)^{p/2}} .\label{BDG}
\end{align} 
\end{proposition}

Let us recall the following theorem  (see \cite{Da}) on the existence of a unique mild solution to equation \eqref{PAM}. 
 
\begin{proposition}  \label{prop2.2}
There exists a unique measurable and  adapted random field $$u=\left\{u(t,x):(t,x)\in  \R_+\times \R^d\right\}$$ such that for all $T>0$ and $p\ge 2$
\begin{equation}  \label{estu1}
 \sup_{(t,x) \in [0,T] \times\R^d }\ex{ |u(t,x) |^p}  =C_{T,p}<\infty \, ,
 \end{equation}
 and  for all  $t\ge 0$ and $x\in \R^d$
\begin{align}\label{mild}
u(t,x)=
1+\int_{[0,t]\times \R^d} \bm{p}_{t-s}(x-y)u(s,y)W(ds,dy).
\end{align}
\end{proposition}

For the existence of densities and the estimates, we will also need some non-degeneracy conditions. For this we will need the following result, which is a corollary  \cite[Theorem 1.5]{ChHuNu21}.
\begin{proposition} \label{prop:negativeu} Let $u$ be the mild solution to the equation \eqref{mild}. Then for all $p>0$, and $t>0$: \begin{align*}
  \sup_{0\le s\le t}  \ex{\left|u(s,0)\right|^{-p}} =\kappa_{t,p}<+\infty.
\end{align*}
\end{proposition}
 
 \subsection{Malliavin Calculus}
   For a smooth and cylindrical random variable of the form $F= f(W(\varphi_1), \dots , W(\varphi_n))$, with $\varphi_i \in \mathfrak{H}$, $1\le i\le n$,  and $f \in C_b^{\infty}(\mathbb{R}^n)$ ($f$ and its partial derivatives are bounded), we define its Malliavin derivative as the $\mathfrak{H}$-valued random variable given by
\[
 DF = \sum_{i=1}^n \frac{\partial f}{\partial x_i} (W(\varphi_1), \dots, W(\varphi_n))\varphi_i\ .
\]
By iteration, we can also define the $k$-th derivative $D^k F$, which is an element in the space $L^2(\Omega; \mathfrak{H}^{\otimes k})$. For any real $p\ge 1$ and any integer $k\ge 1$, the Sobolev space $\mathbb{D}^{k,p}$ is defined as the closure of the space of smooth and cylindrical random variables with respect to the norm $\|\cdot\|_{k,p}$ defined by 
\[
 \|F\|^p_{k,p} = \ex{|F|^p} + \sum_{i=1}^k \ex{\|D^i F\|^p_{\mathfrak{H}^{\otimes i}}}.
\]
We define the divergence operator $\delta$ as the adjoint in $L^2$ of the derivative operator $D$. Namely, an element $U \in L^2(\Omega; \mathfrak{H})$ belongs to the domain of $\delta$, denoted by $\dom \delta$, if there is a constant $c_U > 0$ depending on $U$ satisfying 
\[
|\ex{\langle DF, U\rangle_{\mathfrak{H}}}| \leq c_U \|F\|_{L^2(\Omega)}
\] 
for any $F \in \mathbb{D}^{1,2}$.  If $U\in \dom \,\delta$, the random variable $\delta(U)$ is defined by the duality relationship 
\begin{equation} \label{dua}
\ex{F\delta(U)} = \ex{\langle DF, U \rangle_{\mathfrak{H}}} \, ,
\end{equation}
which is valid for all $F \in \mathbb{D}^{1,2}$. 
 
The It\^o-Walsh stochastic integral we introduced in the previous subsection coincides with the divergence operator for adapted processes (see, for instance, \cite{ChKhNuPu21}). That is to say, for any adapted random field $U= \{ U(t,x): (t,x) \in \mathbb{R}_+\times \mathbb{R}\}$ such that
$\ex{ \| U \|^2_{\mathfrak{H}}}<\infty$, 
  we have that $U\in \dom \delta$ and
   \begin{align}\label{adapted}
    \delta(U)=\int_{\R_+\times\R^d}U(s,y)W(ds,dy).
\end{align} 
The operators $D$ and $\delta$ satisfy the following commutation relation 
\begin{align}
\label{commute} D_{s,y}(\delta(V))=V(s,y)+\delta(D_{s,y} V), \quad
\end{align}
for almost all $(s,y) \in \R_+\times \R^d$, provided
   $V\in \D^{1,2}(\Omega; \HH)$  is such that for almost all $(s,y) \in \R_+\times \R^d$,
   $D_{s,y} V$ belongs to the domain of the divergence in $L^2$
   and 
   \[
   \ex{   \int_{\R_+\times\R^{2d}}  | \delta(D_{s,y} V)| | \delta(D_{s,y'} V)| |y-y'|^{-\beta}dydy'ds } <\infty
   \]
    (see \cite[Proposition 1.3.2]{Nu}).

 The following estimate can be found in \cite{KuNu2022} (see also \cite{ChHuNu21} for the differentiability of the solution). 
\begin{proposition}   \label{prop1}  Let $u$ be the solution to the stochastic heat equation \eqref{PAM}. 
 Fix $(t,x)\in \mathbb{R}_+\times \R^d$. Then $u(t,x) \in  \cap_{p\ge 2} \D^{2,p}$. Moreover,
 \begin{itemize}
\item[(a)] \,The derivative $D_{s,y}u(t,x)$ satisfies \begin{align}\label{Du}
    D_{s,y}u(t,x)=\bm{p}_{t-s}(x-y)u(s,y)+\int_{[s,t]\times \R^d} \bm{p}_{t-\tau}(x-\xi)D_{s,y}u(\tau,\xi) W(d\tau,d\xi)
\end{align} 
if $s<t$ and $D_{s,y}u(t,x)=0$ if $s>t$. For all $0<s<t<T$, we have \begin{align}\label{estDu}
    \left\|D_{s,y}u(t,x) \right\|_p \leq C_{T,p} \bm{p}_{t-s}(x-y),
\end{align}
where  $C_{T,p}$ is a constant that depends  on $T$ and $p$. 
\item[(b)]\, The second derivative $D_{r,z}D_{s,y}u(t,x)$ satisfies
 \begin{align}\label{DDu}
   D_{r,z} D_{s,y}u(t,x) &=\bm{p}_{t-s}(x-y)
   D_{r,z}u(s,y)\\
   &+ \nonumber \int_{[s,t]\times \R^d} \bm{p}_{t-\tau}(x-\xi)D_{r,z}D_{s,y}u(\tau,\xi) W(d\tau,d\xi)
\end{align} 
if $r<s<t$. For all $0<r<s<t<T$, we have
\begin{align}\label{estDDu}
    \left\| D_{r,z}D_{s,y}u(t,x) \right\|_{p} \leq C_{T,p} 
    \bm{p}_{t-s}(x-y) \bm{p}_{s-r}(y-z),
\end{align}
where  $C_{T,p}$ is a constant that depends  on $T$ and $p$. 
\end{itemize}
\end{proposition}

 \subsection{Malliavin-Stein Method}
 In this subsection we recall a result which gives a bound on the uniform distance between the density of a random variable and the density of the standard normal distribution under some assumptions. A version of this result first appeared in  \cite{HuLuNu} and the proof of Theorem  \ref{densityapprox2} below can be found in \cite{KuNu21}. For the existence of the density, see also \cite[Proposition 1]{CaFeNu}.
 
 For $V\in L^2(\Omega ; \mathfrak{H})$ and $F\in \mathbb{D}^{1,2}$ we will make use of the notation $D_VF= \langle DF, V \rangle_{\mathfrak{H}}$.
 \begin{theorem}
\label{densityapprox2} 
Assume that $V \in \D^{1,6}(\Omega; \H)$ and $F=\delta(V) \in \D^{2,6}$ with $\ex{F}=0$, $\ex{F^2}=1$ and $\left( D_VF \right)^{-1} \in L^{4}(\Omega)$. Then, $V/D_VF \in {\rm Dom} \,  \delta$, $F$ admits a density $f_F(x)$ and the following inequality holds true
   \begin{align}  \nonumber
  \sup_{x\in \R}|f_F(x)-\phi(x)| \leq  & \left(\left\| {F}\right\|_{4}\left\| (D_VF)^{-1} \right\|_{4}+2\right) \left\|  1-D_VF\right\|_2\\ 
  \label{e85}
   &+\left\| (D_VF)^{-1} \right\|_{4}^2 \left\|D_V\left({D_VF}\right)\right\|_{2},
\end{align}
where $\phi(x)$ is the density of  the law $ N(0,1)$. 
\end{theorem}

 \subsection{Order of variance and total variation distance}

Let $u$ be the mild  solution to \eqref{PAM} with $u_0=1$ given in Proposition \ref{prop2.2}.
 Then, for any fixed $t>0$,  the  random variable $F_{R,t}$ defined in \eqref{FR} is given by
    \begin{align*} 
   F_{R,t} &=\frac{1}{\sigma_{R,t}}\left(\int_{ \bm{Q}_R} \int_{[0,t] \times \R^d} \bm{p}_{t-\tau}(x-\xi)u(\tau,\xi)W\left(d\tau,d\xi\right)dx\right) .
   \end{align*}
 Using  the stochastic Fubini's theorem, we can rewrite $F_{R,t}$ as  follows 
   \begin{align*}
F_{R,t} &=
\int_{[0,t] \times \R^d}  \frac{1}{\sigma_{R,t}}  \left(\int_{\bm{Q}_R} \bm{p}_{t-\tau}(x-\xi)u(\tau,\xi)dx\right) W\left(d\tau,d\xi\right)\\&=
 \int_{[0,t] \times \R^d}  \varphi_{R,t}(\tau,\xi) u(\tau,\xi)W(d\tau,d\xi),
\end{align*}
where \begin{align}\label{varphi}  \varphi_{R,t}(\tau,\xi) := \frac 1{ \sigma_{R,t}}\int_{\bm{Q}_R} {\bm{p}_{t-\tau}(x-\xi)} dx.
 \end{align}
So, using \eqref{adapted}, we obtain
 \[
 F_{R,t}=\delta(V_{R,t}),
 \]
where
\begin{equation} \label{v}
    V_{R,t}(\tau,\xi)= \mathbf{1}_{[0,t]} (\tau) \varphi_{R,t}(\tau,\xi) u(\tau,\xi).
\end{equation}

 The following result provides the asymptotic behavior of $\sigma^2_{R,t}$ as $R\to \infty$ and an upper bound for the total variation distance between the law of $F_{R,t}$ and the one-dimensional standard distribution (see \cite[Theorem 1.1]{HuNuViZh20}).
 
 \begin{theorem}  Let  $u=\left\{u(t,x): (t,x)\in  \R_+\times \R^d\right\}$ be the mild solution to the stochastic heat equation \eqref{PAM} with the initial condition $u_0=1$. 
Fix $t>0$ and let $F_{R,t}$, and $\sigma_{R,t}$ be defined as in \eqref{FR}. Then, 
\begin{align} \label{variance}
  \lim_{R\rightarrow \infty}   \frac{\sigma^2_{R,t}}{R^{2d-\beta}} = k_\beta \,t,
   \end{align}
    where \begin{align*}
        k_{\beta}:=\int_{\bm{Q}_1^2}|x_1-x_2|^{-\beta}dx_1dx_2.
    \end{align*}
  Moreover, for all $R>0$, 
\begin{align*}
  \d_{\TV}(F_{R,t},Z) \leq \ C_{t,\beta} R^{-\beta/2},
     \end{align*} 
     where Z is a $N(0,1)$ random variable.
\end{theorem}

 \section{Proof of Theorem \ref{thm:rate1}}
In order to apply Theorem~\ref{densityapprox2}, we need estimates on the negative moments of $D_{V_{R,t}}F_{R,t}$. The next proposition provides this type of estimates.

\begin{proposition}\label{F^-1} 
Let $u$ be the mild solution to stochastic the heat equation \eqref{PAM} with $u_0=1$, let $F_{R,t}$ be  defined  as in \eqref{FR} and fix $p\geq 2$. Then, 
\begin{equation}
\sup_{R\ge 1} \ex{\left| D_{V_{R,t}}F_{R,t} \right|^{-p} }<\infty.  \label{b1}
\end{equation}
\end{proposition}
 
 \begin{proof}
 Consider the Malliavin derivative of $F_{R,t}$ given by
  \begin{equation} \label{EQU1}
    D_{s,y}F_{R,t}=\frac{1}{\sigma_{R,t}} \int_{\bm{Q}_R} D_{s,y}u(t,x) dx.
\end{equation}
  From \eqref{v} and \eqref{EQU1}, we can write
 \begin{align}  \nonumber
D_{V_{R,t}}F_{R,t}&= \int_0^t \int_{\R^{2d}}  V_{R,t}(s,\ytilde)D_{s,y}F_{R,t} |y-\ytilde|^{-\beta} d\ytilde dyds
 \\&=\frac{1}{\sigma_{R,t}^2}\int_{\bm{Q}_R^2} \int_0^t \int_{\R^{2d}} \bm{p}_{t-s}(x_1-\ytilde)u(s,\ytilde)D_{s,y}u(t,x_2)|y-\ytilde|^{-\beta} d\ytilde\,dydsdx_1 dx_2.  \label{posit}
\end{align}
First note that the integrand in \eqref{posit} is nonnegative since for all $(t,x)\in \R_+\times\R^d$, $u(t,x)\geq 0$, and for all $(s,y),(t,x)\in \R_+\times\R^d$, $D_{s,y}u(t,x) \geq 0$ a.s. 


Fix $t>\epsilon>0$ and set $\talpha:=t-\e^{\alpha}$ where $\alpha<1$.  Then we can estimate \eqref{posit} as follows: 
\begin{align}  \nonumber
D_{V_{R,t}}F_{R,t}&\geq \frac{1}{\sigma_{R,t}^2}\int_{\bm{Q}_R^2}\int_{\talpha}^t\int_{\R^{2d}}\bm{p}_{t-s}(x_1-\ytilde)u(s,\ytilde)D_{s,y}u(t,x_2)|y-\ytilde|^{-\beta} d\ytilde dydsdx_1 dx_2
 \\&=:\mathrm{\bm{I}}_{R,\e}.\nonumber
\end{align}
 Using this estimate, we get
\[
\pr{  D_{V_{R,t}}F_{R,t} <\e } \leq   \pr{ \mathrm{\bm{I}}_{R,\e} < \e }. 
\]
From \eqref{Du}, we obtain
\begin{align*}
\mathrm{\bm{I}}_{R,\e}  &=  \int_{\talpha}^t\int_{\R^{2d}} \varphi_{R,t}(s,y)\varphi_{R,t}(s,\ytilde)u(s,z)u(s,\ytilde)|y-\ytilde|^{-\beta} d\ytilde dyds\\ 
 &+\int_{\talpha}^t \int_{\R^{2d}} \left(\int_{[s,t]\times\R^d}\varphi_{R,t}(\tau,\xi)  D_{s,y}u(\tau,\xi)W(d\tau,d\xi)\right)\varphi_{R,t}(s,\ytilde) u(s,\ytilde)|y-\ytilde|^{-\beta} d\ytilde dyds\\
 &=:\mathrm{\bm{I}}_1+\mathrm{\bm{I}}_2.
\end{align*}
Taking into account the estimates
\begin{align}
  \label{prob} \pr{\mathrm{\bm{I}}_1+\mathrm{\bm{I}}_2 < \e}&\leq   \pr{\mathrm{\bm{I}}_1 < 2\e}+\pr{\mathrm{\bm{I}_1}+\mathrm{\bm{I}}_2 < \e, \, \mathrm{\bm{I}}_1\geq 2\e}\\ \nonumber
    &\leq \pr{\mathrm{\bm{I}}_1 <2 \e}+\pr{ |\mathrm{\bm{I}}_2|> \e},  
\end{align}
 we have
\begin{align*}
 \pr{  \mathrm{\bm{I}}_{R,\e} < \e } \leq  \pr{\mathrm{\bm{I}}_1 <2 \e} +\pr{\left|  \mathrm{\bm{I}}_2 \right|> \e  }.
\end{align*}
We shall next estimate these probabilities, starting with the first one.
\begin{align*}
 K:=\pr{ \mathrm{\bm{I}}_1< 2\e} &=  \pr{ \int_{\talpha}^t\int_{\R^{2d}} \varphi_{R,t}(s,y)\varphi_{R,t}(s,\ytilde)u(s,y)u(s,\ytilde)|y-\ytilde|^{-\beta} d\ytilde dyds < 2\e }.
\end{align*}
By Chebyshev's inequality, for any $q\geq 2$ we obtain
\begin{align}  \nonumber
 K & =\P\left(   \left[ \int_{\talpha}^t\int_{\R^{2d}} \varphi_{R,t}(s,y)\varphi_{R,t}(s,\ytilde)u(s,y)u(s,\ytilde)|y-\ytilde|^{-\beta} d\ytilde dyds\right] ^{-1} > \frac 1{2\e} \right) \\
 &  \le (2\e)^q  \ex{ \left( \int_{\talpha}^t\int_{\R^{2d}} \varphi_{R,t}(s,y)\varphi_{R,t}(s,\ytilde)u(s,y)u(s,\ytilde)|y-\ytilde|^{-\beta} d\ytilde dyds \right) ^{-q}}.  \label{e8}
\end{align}
Set
\begin{align}\label{m}
m(t_{\alpha}, R):=  \int_{\talpha} ^t \int_{\R^{2d}}    \varphi_{R,t}(s,y) \varphi_{R,t}(s,\ytilde) |y-\ytilde|^{-\beta} d\ytilde dyds .
\end{align}
Then, taking into account that the function $x\to x^{-q}$ is convex   and applying Jensen's inequality, we can write
\begin{align}  \nonumber
&  \ex{  \left(\int_{\talpha}^t\int_{\R^{2d}} \varphi_{R,t}(s,y)\varphi_{R,t}(s,\ytilde)u(s,y)u(s,\ytilde)|y-\ytilde|^{-\beta} d\ytilde dyds\right) ^{-q}} \\  \nonumber
&= m(t_{\alpha}, R)^{-q}  \ex{ \left(  \frac 1{ m(t_{\alpha},R)}  \int_{\talpha}^t\int_{\R^{2d}} \varphi_{R,t}(s,y)\varphi_{R,t}(s,\ytilde)u(s,y)u(s,\ytilde)|y-\ytilde|^{-\beta} d\ytilde dyds \right)^{-q}} \\  \label{jensen}
& \le   {m(t_{\alpha}, R)^{-q-1} }   \int_{\talpha}^t\int_{\R^{2d}} \varphi_{R,t}(s,y)\varphi_{R,t}(s,\ytilde)\ex{\left(u(s,y)u(s,\ytilde)\right)^{-q}}|y-\ytilde|^{-\beta} d\ytilde dyds.
\end{align}
Using H\"older's inequality, the  stationarity of solution in  the space variable (see, for insance, \cite[Lemma 7.1]{ChKhNuPu21}), and finally Proposition~\ref{prop:negativeu}, we get
 \begin{align}
 \label{uq}\ex{\left(u(s,y)u(s,\ytilde)\right)^{-q}} \leq \sqrt{\ex{u(s,y)^{-2q}}\ex{u(s,\ytilde)^{-2q}}}   =\ex{u(s,0)^{-2q}} \leq \kappa_{t,2q}
\end{align}
for $t_{\alpha} \leq s\leq t$. Substituting the estimate \eqref{uq} into \eqref{jensen}, we get
 \begin{align} \nonumber
   & \ex{  \left(\int_{\talpha}^t\int_{\R^{2d}} \varphi_{R,t}(s,y)\varphi_{R,t}(s,\ytilde)u(s,y)u(s,\ytilde)|y-\ytilde|^{-\beta} d\ytilde dyds\right) ^{-q}} \\  &\le \kappa_{t,2q}{m(t_{\alpha}, R)^{-q-1} }   \int_{\talpha}^t\int_{\R^{2d}} \varphi_{R,t}(s,y)\varphi_{R,t}(s,\ytilde)|y-\ytilde|^{-\beta} d\ytilde dyds
=\kappa_{t,2q}{m(t_{\alpha}, R)^{-q} }.  \label{EQU2}
\end{align}
Therefore, from \eqref{e8} and \eqref{EQU2}, we obtain
\begin{equation} \nonumber
K \le  \kappa_{t,q2}(2\e)^q  m(t_{\alpha}, R)^{-q}.
\end{equation}
Using Lemma \ref{lem:m}, we obtain \begin{equation} \label{estK}
K \le  C_{t,q, \beta, d}\epsilon^{(1-\alpha)q},
\end{equation}
where $C_{t,q, \beta, d}$ is a constant depending on $t$, $q$, $\beta$ and $d$.
Next, we will estimate the probability
$L:=\P\left(\left|\mathrm{\bm{I}}_2\right|> \e \right)$. 
Using the stochastic Fubini's theorem and Chebyschev's inequality, for any $q\geq 2$ we have 
\begin{align*}
   &L \leq \frac{1}{\e^q}\\
   & \times \ex{\left|\int_{\talpha}^t \int_{\R^{d}}\left(\int_{\R^{2d}}\int_{\talpha}^{\tau}  \varphi_{R,t}(s,\ytilde)\varphi_{R,t}(\tau,\xi)D_{s,y}u(\tau,\xi)u(s,\ytilde)|y-\ytilde|^{-\beta}dsdyd\ytilde\right) W(d\tau,d\xi)\right|^q}.
\end{align*}
Then, applying  Burkholder-Davis-Gundy inequality \eqref{BDG}, followed by Minkowski's inequality,  we get
\begin{align}   \nonumber
  & L   \leq \frac{C_q}{\e^q}\E\left[\left|\int_{\talpha}^t \int_{\R^{6d}}\int_{\talpha}^\tau \int_{\talpha}^\tau  \varphi_{R,t}(s_1,\ytilde_1)\varphi_{R,t}(\tau,\xi) \varphi_{R,t}(s_2,\ytilde_2)\varphi_{R,t}(\tau,\xitilde) D_{s_1,y_1}u(\tau,\xi)u(s_1,\ytilde_1)  \right.\right.\\ & \nonumber \hskip 1cm
  \times \left.\left.D_{s_2,y_2}u(\tau,\xitilde)u(s_2,\ytilde_2) |\xi-\xitilde|^{-\beta} |y_1-\ytilde_1|^{-\beta} |y_2-\ytilde_2|^{-\beta} ds_1ds_2 d\ytilde_1d\ytilde_2dy_1dy_2 d\xitilde d\xi  d\tau \right|^{\frac q2}\right]    
    \\ \nonumber &=\frac{C_q}{\e^q} \E \Bigg[ \Bigg|\int_{\talpha}^t\int_{\talpha}^\tau \int_{\talpha}^\tau\int_{\R^{6d}}  \varphi_{R,t}(s_1,\ytilde_1)\varphi_{R,t}(\tau,\xi) \varphi_{R,t}(s_2,\ytilde_2)\varphi_{R,t}(\tau,\xitilde)  X_{s_1,y_1,\ytilde_1,s_2,y_2,\ytilde_2}(\tau,\xi,\xitilde)\\  \nonumber
    & \qquad  \qquad  \qquad  \qquad \qquad \times |\xi-\xitilde|^{-\beta} |y_1-\ytilde_1|^{-\beta} |y_2-\ytilde_2|^{-\beta} d\ytilde_1 d\ytilde_2dy_1dy_2 d\xitilde d\xi  ds_1ds_2 d\tau    \Bigg|^{\frac q2} \Bigg]
    \\ \nonumber
     &    \leq \frac{C_q}{\e^q}\Bigg( \int_{\talpha}^t\int_{\talpha}^\tau \int_{\talpha}^\tau\int_{\R^{6d}} \varphi_{R,t}(s_1,\ytilde_1)\varphi_{R,t}(\tau,\xi) \varphi_{R,t}(s_2,\ytilde_2)\varphi_{R,t}(\tau,\xitilde)\left\|X_{s_1,y_1,\ytilde_1,s_2,y_2,\ytilde_2}(\tau,\xi,\xitilde) \right\|_{q/2} \\
    & \qquad  \qquad  \qquad  \qquad \qquad    \times  |\xi-\xitilde|^{-\beta} |y_1-\ytilde_1|^{-\beta} |y_2-\ytilde_2|^{-\beta} d\ytilde_1d\ytilde_2dy_1dy_2 d\xitilde d\xi 
     ds_1ds_2 d\tau \Bigg)^{\frac q2}, \label{BDG1}
\end{align}
where 
\begin{align*}
X_{s_1,y_1,\ytilde_1,s_2,y_2,\ytilde_2}(\tau,\xi,\xitilde):=& D_{s_1,y_1}u(\tau,\xi)u(s_1,\ytilde_1)D_{s_2,y_2}u(\tau,\xitilde)u(s_2,\ytilde_2) .
\end{align*}
Using H\"older's inequality,  the estimate \eqref{estDu}  and the fact that $\sup_{(r,z) \in [0,t]\times\R}\|u(r,z)\|_{q} <\infty$ for all
$q\ge 2$,   we have
 \[
  \|X_{s_1,y_1,\ytilde_1,s_2,y_2,\ytilde_2}(\tau,\xi,\xitilde) \|_{q/2} \leq C \bm{p}_{\tau-s_1}(y_1-\xi) \bm{p}_{\tau-s_2}(y_2-\xitilde),
\]
where here and in the rest of the proof $C$ will denote a generic constant that depends on $t$, $q$, $\beta$ and $d$.
Plugging this bound  in the estimate \eqref{BDG1}, we see that
\begin{align}  \nonumber
  L& \leq   \frac{C}{\e^q} \Big( \int_{\talpha}^t\int_{\talpha}^{\tau} \int_{\talpha}^{\tau}\int_{\R^{6d}}  \bm{p}_{\tau-s_1}(y_1-\xi) \bm{p}_{\tau-s_2}(y_2-\xitilde) \varphi_{R,t}(s_1,\ytilde_1)\varphi_{R,t}(\tau,\xi) \varphi_{R,t}(s_2,\ytilde_2)\varphi_{R,t}(\tau,\xitilde)\\
    & \qquad  \qquad  \qquad  \qquad \times   |\xi-\xitilde|^{-\beta} |y_1-\ytilde_1|^{-\beta} |y_2-\ytilde_2|^{-\beta}d\ytilde_1d\ytilde_2dy_1dy_2 d\xitilde d\xi  ds_1ds_2 d\tau
    \Big)^{\frac q2} . \label{L}
\end{align}
 Then, we can write
\begin{align*} 
 L& \leq   \frac{C}{\e^q} \left(\frac{1}{\sigma^4_{R,t}}  \int_{t_{\alpha}}^t \int_{t_{\alpha}}^{\tau}\int_{t_{\alpha}}^{\tau}  E_{R,t}(s_1,s_2,\tau)ds_1ds_2d\tau  \right)^{\frac q2},
\end{align*}
where 
\begin{align} \label{ert}
    E_{R,t}(s_1,s_2,\tau)\nonumber 
    &   :=\int_{\bm{Q}_R^4}dx_1dx_2dx_3dx_4 \int_{\R^{6d}} d\ytilde_1d\ytilde_2 d\xi  \xitilde dy_1dy_2
      |\xi-\xitilde|^{-\beta} |y_1-\ytilde_1|^{-\beta}|y_2-\ytilde_2|^{-\beta} \\ \nonumber
      & \qquad \times \, \bm{p}_{t-\tau}(x_1-\xi)\bm{p}_{t-\tau}(x_2-\xitilde) \bm{p}_{t-s_1}(x_3-\ytilde_1) \\
      & \qquad \times  \bm{p}_{t-s_2}(x_4-\ytilde_2)\bm{p}_{\tau-s_1}(\xi-y_1) \bm{p}_{\tau-s_2}(y_2-\xitilde).
\end{align}
Now, using the estimate in Lemma~\ref{lem:E} and \eqref{variance}, we get, for $R\ge 1$,
\begin{align}\label{estL}
   L \leq \frac{C}{\e^q} \left(\frac{R^{4d-3\beta}}{\sigma_{R,t}^4}\int_{t_{\alpha}}^t \int_{t_{\alpha}}^{\tau}\int_{t_{\alpha}}^{\tau} ds_1ds_2d\tau\right)^{q/2} \leq C_{}\e^{(\frac{3\alpha}{2}-1)q}.
\end{align}
Finally,  choosing  $\alpha=4/5$ in \eqref{estK} and \eqref{estL}, we obtain
\begin{align*}
     \sup_{R\ge 1} \P\left(D_{V_{R,t}}F_{R,t} < \e \right) \leq C_{} \e^{q/5}.
\end{align*} 
Together with Lemma~\ref{lem:conditionnegativemoments}, this estimate completes our proof.
\end{proof}

Now, we are ready to prove the main result.

\begin{proof}[Proof of Theorem \ref{thm:rate1}]
We will apply Theorem \ref{densityapprox2} to the random variable  $F_{R,t} = \delta (V_{R,t})$. 
Fix  $t>0$. 
Along the proof $C$ will denote a generic constant that depends on $t$, $q$, $\beta$ and $d$.
It has been already proven in \cite{HuNuViZh20} that
\begin{equation} \label{dv}
\left\|1-D_{V_{R,t}}F_{R,t}\right\|_{2}\leq C\,{R^{-\beta/2}}.
\end{equation}
Taking into account Proposition~\ref{F^-1} we are only left to estimate the term $\left\| D_{V_{R,t}}\left({D_{V_{R,t}}F_{R,t}}\right)\right\|_{2}$
appearing in \eqref{e85}.  Recall that

\begin{align*}
   D_{V_{R,t}}F_{R,t}&= \frac{1}{\sigma_{R,t}}\int_0^t\int_{\R^{2d}}\int_{\bm{Q}_R} \varphi_{R,t}(s,\ytilde)u(s,\ytilde) D_{s,y} u(t,x)  |y-\ytilde|^{-\beta} dx d\ytilde dy ds.
\end{align*}
Taking  the Malliavin derivative,  we get 
\begin{align*}
D_{r,z} \left(D_{V_{R,t}}F_{R,t}\right)=\frac{1}{\sigma_{R,t} }\int_{r}^t\int_{\R^{2d}}\int_{\bm{Q}_R} \varphi_{R,t}(s,\ytilde) D_{r,z} u(s,\ytilde) D_{s,y} u(t,x)|y-\ytilde|^{-\beta} dx d\ytilde dy ds
\\ +\frac{1}{\sigma_{R,t}}\int_0^t\int_{\R^{2d}}\int_{\bm{Q}_R} \varphi_{R,t}(s,\ytilde)u(s,\ytilde) D_{r,z} D_{s,y} u(t,x)|y-\ytilde|^{-\beta} dx d\ytilde dy ds,
\end{align*}
and hence, using \eqref{v}, we can write
\begin{align*}
D_{V_{R,t}} \left(D_{V_{R,t}}F_{R,t}\right)&=\frac{1}{\sigma_{R,t}}\int_{0}^t\int_{r}^t\int_{\R^{4d}}\int_{\bm{Q}_R} \varphi_{R,t}(r,\ztilde)\varphi_{R,t}(s,\ytilde)u(r,\ztilde) D_{r,z} u(s,\ytilde) D_{s,y} u(t,x)\\
& \qquad  \times |y-\ytilde|^{-\beta} |z-\ztilde|^{-\beta}  dx d\ytilde d\ztilde dz dy dsdr
\\ & \qquad +\frac{2}{\sigma_{R,t}}\int_{0}^t\int_{r}^{t}\int_{\R^{4d}}\int_{\bm{Q}_R}\varphi_{R,t}(r,\ztilde)\varphi_{R,t}(s,\ytilde)u(r,\ztilde)u(s,\ytilde) D_{r,z} D_{s,y} u(t,x)\\
&  \qquad \times |y-\ytilde|^{-\beta} |z-\ztilde|^{-\beta}  dx d\ytilde d\ztilde dz dy dsdr.
\end{align*}
Now using \eqref{Du} for $D_{s,y}u(t,x)$ and \eqref{DDu}  for $D_{r,z}D_{s,y}u(t,x)$, respectively, we have 
\begin{align*}
D_{V_{R,t}} \left(D_{V_{R,t}}F_{R,t}\right)= 
\y^0_{R,t}+\y^1_{R,t}+\y^2_{R,t},
\end{align*}
where: \begin{align*}
    \y^0_{R,t}&=2 \int_{0}^t\int_{r}^t\int_{\R^{4d}} \varphi_{R,t}(r,\ztilde)\varphi_{R,t}(s,\ytilde)\varphi_{R,t}(s,y) D_{r,z} u(s,\ytilde) u(r,\ztilde)u(s,y)\\
    & \qquad \times |y-\ytilde|^{-\beta} |z-\ztilde|^{-\beta}   d\ytilde d\ztilde dz dy dsdr , \\
    & \\
    \y^1_{R,t}&= \int_{0}^t\int_{r}^t\int_{\R^{4d}}  \left(\int_{[s,t]\times\R^d}\varphi_{R,t}(\tau,\xi) D_{s,y}u(\tau,\xi) W(d\tau,d\xi) \right)
    \varphi_{R,t}(r,\ztilde)\varphi_{R,t}(s,\ytilde)\\
    & \qquad  \times u(r,\ztilde) D_{r,z} u(s,\ytilde)
     |y-\ytilde|^{-\beta} |z-\ztilde|^{-\beta}   d\ytilde d\ztilde dz dy dsdr,   \\
     & \\
    \y^2_{R,t}&=2\int_{0}^t\int_r^{t}\int_{\R^{4d}}  \left(\int_{[ s, t]\times\R^d}\varphi_{R,t}(\tau,\xi) D_{r,z}D_{s,y}u(\tau,\xi) W(d\tau,d\xi) \right)
     \varphi_{R,t}(r,\ztilde)\varphi_{R,t}(s,\ytilde)\\ 
     & \qquad \times u(r,\ztilde)u(s,\ytilde)
    |y-\ytilde|^{-\beta} |z-\ztilde|^{-\beta}   d\ytilde d\ztilde dz dy dsdr.
\end{align*}
We will handle the terms $\y^1_{R,t}$ and $\y^2_{R,t}$ simultaneously. To this end, for $i=1,2$, let
\begin{align*}
 \y^i_{R,t} = 
   \int_{[0,t]\times\R^d}    
  &\Bigg( \int_0^{\tau}  \int_{r} ^{\tau} \int _{\R^{4d}}    \varphi_{R,t}(r,\ztilde)\varphi_{R,t}(s,\ytilde) Z^i_{r,z,\ztilde, s,y,\ytilde}(\tau ,\xi ) \\
  & \qquad \times  |y-\ytilde|^{-\beta} |z-\ztilde|^{-\beta} d\ytilde d\ztilde dz dy dsdr  \Bigg)  \varphi_{R,t} (\tau,\xi) W(d\tau,d\xi),
\end{align*}
where 
\begin{align} \label{Z} 
Z^i_{r,z,\ztilde, s,y,\ytilde}(\tau ,\xi )  &:=  \begin{cases} &u(r,\ztilde) D_{r,z} u(s,\ytilde)D_{s,y}u(\tau,\xi),  \text{ for }i=1  \\  
& u(r,\ztilde)u(s,y)   D_{r,z} D_{s,y} u(\tau,\xi), \text{ for }i=2 .
\end{cases}
\end{align}

 \noindent
{\it Estimation of}   $\left\|\y^0_{R,t}\right\|_2$: Note that using the estimates \eqref{estu1} and  \eqref{estDu} and H\"older's inequality we have, for $r < s$,
\[
    \left\| u(r,\ztilde)u(s,y) D_{r,z} u(s,\ytilde) \right\|_2   \leq  C \bm{p}_{s-r}(z-\ytilde).
\]
    As a consequence,
\[
     \left\|\y^0_{R,t}\right\|_2\leq  C \int_{0}^t\int_{r}^t\int_{\R^{4d}} \varphi_{R,t}(r,\ztilde)\varphi_{R,t}(s,\ytilde)\varphi_{R,t}(s,y)\bm{p}_{s-r}(z-\ytilde) |y-\ytilde|^{-\beta} |z-\ztilde|^{-\beta} d\ytilde d\ztilde dz dy dsdr.
\]
    Integrating in the variables $z,y$ and using Lemma \ref{lem:p} and Lemma  \ref{LemmaA.2}, we have
     \begin{align*}
   & \int_{\R^{2d}}\varphi_{R,t}(s,y) \bm{p}_{s-r}(z-\ytilde)|z-\ztilde|^{-\beta}|y-\ytilde|^{-\beta} dz dy\\ &= \frac{1}{\sigma_{R,t}} \int_{\bm{Q}_R}\int_{\R^{2d}} \bm{p}_{t-s}(x-y)\bm{p}_{s-r}(z-\ytilde)|y-\ytilde|^{-\beta}|z-\ztilde|^{-\beta} dz dy dx \\ &\leq 
    C \frac{|\ytilde-\ztilde|^{-\beta}}{\sigma_{R,t}} \int_{\bm{Q}_R}\int_{\R^d}\bm{p}_{t-s}(x-y)|y-\ytilde|^{-\beta}dy dx \leq   C\frac{|\ytilde-\ztilde|^{-\beta}}{\sigma_{R,t}} \int_{\bm{Q}_R} | x-\ytilde|^{-\beta}dx\\
    &\leq C |\ytilde-\ztilde|^{-\beta} \frac{R^{d-\beta}}{\sigma_{R,t}}.
    \end{align*}
 Using this estimate followed with Lemma \ref{lem:varvar} and \eqref{variance}, we get, for $R\ge 1$, 
   \begin{align}\nonumber \label{y0}
     \left\|\y^0_{R,t}\right\|_2 &\leq \frac{ C R^{d-\beta}}{\sigma_{R,t}}\int_0^t \int_{r}^t  \int_{\R^{2d}} \varphi_{R,t}(r,\ztilde)\varphi_{R,t}(s,\ytilde)|\ztilde-\ytilde|^{-\beta}d\ytilde d\ztilde  dsdr
  \\ & \leq  \frac{ R^{d-\beta} C}{\sigma_{R,t}} \leq CR^{-\beta/2}.
\end{align}

      \medskip
 \noindent
{\it Estimation of}   $\left\|\y^i_{R,t}\right\|_2$ for $i=1,2$:  Using  the It\^o-Walsh isometry of the stochastic integral and Cauchy-Schwarz inequality, we obtain
 \begin{align*}
       \left\|\y^i_{R,t}\right\|_2^2   &= \int_0^t \int_{\R^{2d}}    \ex{ I(Z^i(\tau,\xi)) I(Z^i(\tau,\tilde{\xi}))}  \varphi_{R,t} (\tau,\xi)\varphi_{R,t} (\tau,\tilde{\xi})  |\xi-\tilde{\xi}|^{-\beta}d\xi d\xitilde d\tau \\
       &\leq \int_{0}^t \int_{\R^{2d}} \|I(Z^{i}(\tau,\xi))\|_2\|I(Z^i(\tau,\tilde{\xi}))\|_2 \varphi_{R,t} (\tau,\xi)\varphi_{R,t} (\tau,\tilde{\xi})|\xi-\tilde{\xi}|^{-\beta}d\xi d\xitilde d\tau
  \end{align*}
  where
   \begin{align*}
  I(Z^i(\tau,\xi))= \int_0^\tau \int _r^\tau \int_{\R^{4d}}   \varphi_{R,t}(s,\ytilde)\varphi_{R,t}(r,\ztilde)    Z^i_{r,z,\ztilde, s, y,\ytilde}(\tau,\xi)|y-\ytilde|^{-\beta} |z-\ztilde|^{-\beta}dy dz d\ytilde  d\ztilde ds dr.
  \end{align*}
From the definition \eqref{Z}, using H\"older's inequality and the estimates \eqref{estu1}, \eqref{estDu} and \eqref{estDDu}, we have 
\[
    \|Z^i_{r,z,\ztilde, s, y,\ytilde}(\tau,\xi) \|_2 \leq  C\phi^{(i)}_{r,z,\ztilde, s, y,\ytilde}(\tau,\xi),
    \]
    where
\begin{align} \label{phi}
    \phi^{(i)}_{r,z,\ztilde, s, y,\ytilde}(\tau,\xi)= \begin{cases}
    \bm{p}_{s-r} (\ytilde-z) \bm{p}_{\tau-s} ( \xi-y), &\text{ for }i=1 \\ \bm{p}_{\tau-s} (\xi-y) \bm{p}_{s-r} ( y-z), &\text{ for }i=2.
    \end{cases}
\end{align}
Hence, we get
\begin{align*}
     &  \left\|\y^i_{R,t}\right\|_2^2    \le   C\int_0^t \int_{\R^{2d}}     \varphi_{R,t} (\tau,\xi) \varphi_{R,t} (\tau,\tilde{\xi})  \Phi^{(i)}(\tau, \xi) \Phi^{(i)}(\tau, \tilde{\xi}) |\xi-\tilde{\xi}|^{-\beta}d\xi d\xitilde d\tau,
       \end{align*}
     where 
     \begin{equation} \label{Phi}
   \Phi^{(i)}(\tau, \xi)=\int_0^\tau \int _r^\tau\int_{\R^{4d}}   \varphi_{R,t}(s,\ytilde)\varphi_{R,t}(r,\ztilde) \phi^{(i)}_{r,z,\ztilde, s, y,\ytilde}(\tau,\xi)|z-\ztilde|^{-\beta} |y-\ytilde|^{-\beta}d\ytilde d\ztilde   ds dr.
       \end{equation}
Using Lemma $\ref{lem:Phib}$, we obtain
\[
       \left\|\y^i_{R,t}\right\|_2^2    \le   C R^{-2\beta}        \int_0^t \int_{\R^{2d}} \varphi_{R,t} (\tau,\xi) \varphi_{R,t} (\tau,\tilde{\xi})    
     \left|\xi-\tilde{\xi}\right|^{-\beta} d\xitilde d\xi d\tau.
\]
Then,  Lemma \ref{lem:varvar} implies that
    \begin{equation} \label{yi}
       \left\|\y^{i}_{R,t}\right\|_2    \le  C R^{-\beta}. 
       \end{equation}
Using the estimates \eqref{y0} and \eqref{yi}, we obtain \begin{align}
    \label{dvdv} \left\|D_{V_{R,t}}\left(D_{V_{R,t}}F_{R,t}\right)\right\|_2\leq C  R^{-\beta/2}.
\end{align}
Finally, plugging \eqref{b1}, \eqref{dv} and \eqref{dvdv} into \eqref{e85} we complete the proof.
 \end{proof}


\appendix 

\section{}

\begin{lemma} \label{lem:p} For any $\beta\in (0,d)$ we have
\begin{align*}
    \sup_{t>0} \int_{\R^d} \bm{p}_{t}(x-y)|y|^{-\beta} dy \leq C|x|^{-\beta},
\end{align*}
for some constant $C>0$ that depends on $d$ and $\beta$.
\end{lemma}
      \begin{proof}  See, for instance,  \cite[Lemma 3.1]{HuNuViZh20}.
      \end{proof}
      
   \begin{lemma} \label{LemmaA.2}
   There is a constant $_{\beta,d}$ depending on $\beta$ and $d$ such that for any $R>0$ we have
   \[
   \int_{Q_R} | x-y|^{-\beta} dy \le C_{\beta,d} R^{d-\beta}.
   \]
   \end{lemma}
    \begin{proof}
    Making the change of variables $y=Rz$, yields
    \[
      \int_{Q_R} | x-y|^{-\beta} dy=R^{d-\beta}   \int_{Q_1} | \frac xR-z|^{-\beta} dz.
      \]
    Then the desired result follows from the fact that 
    \[
    \sup_{w\in \mathbb{R}^d}   \int_{Q_1} | z+w|^{-\beta} dz <\infty.
    \]    
    \end{proof}
    
\begin{lemma}\label{lem:Phib} Let $\Phi^{(i)}$   be as in \eqref{Phi}. There is a constant $C_{t,\beta,d}$ depending on $t$, $\beta$ and $t$, such that for  $0<\tau<t$ and $\xi \in\R^d$ we have
\begin{align*}
     \Phi^i(\tau,\xi) \leq C_{t,\beta,d} R^{-\beta}.
\end{align*}
      \end{lemma}
 
\begin{proof}
We will first consider the case $i=1$. Integrating in $z$ using semigroup property   and then integrating in the variables $y$ and  $\ztilde$, applying Lemma \ref{lem:p},  we have 
\begin{align*}
  &\int_{\R^{4d}}  \varphi_{R,t}(s,\ytilde)\varphi_{R,t}(r,\ztilde) \phi^{(1)}_{r,z,\ztilde, s, y,\ytilde}(\tau,\xi) |y-\ytilde|^{-\beta} |z-\ztilde|^{-\beta} dydz d\ytilde d\ztilde\\
   &=\frac{1}{\sigma_{R,t}^2}\int_{\R^{4d}}\int_{\bm{Q}_R^2}    \bm{p}_{t-s}(x_1-\ytilde)\bm{p}_{t-r}(x_2-z+\ztilde ) \bm{p}_{s-r}(\ytilde-z)  \bm{p}_{\tau-s}(\xi-y) 
   \\
   &\qquad \times |y-\ytilde|^{-\beta} |\ztilde|^{-\beta}dx_1dx_2 dydzd\ytilde d\ztilde 
   \\&=
   \frac{1}{\sigma_{R,t}^2}\int_{\R^{2d}}\int_{\bm{Q}_R^2}    \bm{p}_{t-s}(x_1-\ytilde)\bm{p}_{t+s-2r}(x_2-\ytilde+\ztilde)p_{\tau-s}(\xi-y) |y-\ytilde|^{-\beta} |\ztilde|^{-\beta} d\ytilde d\ztilde
   \\&\leq
\frac{C_{\beta,d}}{\sigma_{R,t}^2}\int_{\R^{d}}\int_{\bm{Q}_R^2}    \bm{p}_{t-s}(x_1-\ytilde)\left|x_2-\ytilde \right|^{-\beta} |\xi-\ytilde|^{-\beta}d\ytilde dx_1dx_2,  
\end{align*}
where $C_{\beta,d}$ is a generic constant depending on $\beta$ and $d$.
Next we apply Lemma \ref{LemmaA.2} to get
\begin{align*}
  &\int_{\R^{4d}}  \varphi_{R,t}(s,\ytilde)\varphi_{R,t}(r,\ztilde) \phi^{(1)}_{r,z,\ztilde, s, y,\ytilde}(\tau,\xi) |y-\ytilde|^{-\beta} |z-\ztilde|^{-\beta} dydz d\ytilde d\ztilde
   \\&\leq
\frac{R^{d-\beta}C_{\beta,d}}{\sigma_{R,t}^2}\int_{\R^{d}}\int_{\bm{Q}_R}    \bm{p}_{t-s}(x_1-\ytilde) |\xi-\ytilde|^{-\beta}d\ytilde dx_1.
\end{align*}
Integrating in the variable $\ytilde$, using Lemma \ref{lem:p} and Lemma \ref{LemmaA.2} and \eqref{variance}, we finally obtain for $0<\tau<t$
    \begin{align*}
     \Phi^i(\tau,\xi) \leq \frac{C_{\beta,d } R^{2d-2\beta}}{\sigma_{R,t}^2} \int_{0}^{\tau}\int_{0}^sdrds \leq  C_{t,\beta,d} R^{-\beta}. 
\end{align*}

Similarly for $i=2$, integrating in $\ytilde, \ztilde$ using Lemma \ref{lem:p}, and then integrating in $x_1,x_2$ using Lemma \ref{LemmaA.2}, we get 
 \begin{align*}
  &\int_{\R^{4d}}  \varphi_{R,t}(s,\ytilde)\varphi_{R,t}(r,\ztilde) \phi^{(2)}_{r,z,\ztilde, s, y,\ytilde}(\tau,\xi) |y-\ytilde|^{-\beta} |z-\ztilde|^{-\beta} dydz d\ytilde d\ztilde\\
   &=\frac{1}{\sigma_{R,t}^2}\int_{\bm{Q}_R^2}   \int_{\R^{4d}} \bm{p}_{t-s}(x_1-\ytilde)\bm{p}_{t-r}(x_2-\ztilde ) \bm{p}_{s-r}(y-z)  \bm{p}_{\tau-s}(\xi-y) |y-\ytilde|^{-\beta} |z-\ztilde|^{-\beta}d\ytilde d\ztilde  dydz dx_1dx_2
   \\&\leq 
   \frac{C_{\beta,d}}{\sigma_{R,t}^2}\int_{\R^{2d}}\int_{\bm{Q}_R^2}    \bm{p}_{s-r}(y-z)\bm{p}_{\tau-s}(\xi-y) |x_1-y|^{-\beta} |x_2-\ztilde|^{-\beta} dx_1dx_2dy dz 
   \\&\leq
\frac{R^{2d-2\beta}C_{\beta,d}}{\sigma_{R,t}^2}\int_{\R^{2d}} \bm{p}_{s-r}(y-z)\bm{p}_{\tau-s}(\xi-y)  dydz   =\frac{R^{2d-2\beta}C_{\beta,d}}{\sigma_{R,t}^2}.
\end{align*}
The desired result follows from \eqref{variance}  after integrating in time variables.

\end{proof}


\begin{lemma}\label{lem:varvar} For $\beta \in (0,d)$ and $t>0$, there is a constant $C_{t,\beta,d}$ depending  on $t$, $\beta$ and $d$ such that  for any $t>0$,
\[
\sup_{s,r \in [0,t]}\int_{\R^{2d}} \varphi_{R,t}(s,y) \varphi_{R,t}(r,z) |y-z| ^{-\beta} dydz  \le C_{t,\beta,d}.
\]
\end{lemma}

\begin{proof} We can write
\begin{align*}
& \int_{\R^{2d}} \varphi_{R,t}(s,y) \varphi_{R,t}(r,z) |y-z| ^{-\beta} dydz  \\
&=\frac 1 {\sigma_{R,t}^2} \int_{Q_R^2} \int_{\R^{2d}} \bm{p}_{t-s}(x_1-y)\bm{p}_{t-r}(x_2-z) |y-z| ^{-\beta} dydz dx_1dx_2
\end{align*} 
Integrating in the variables $y$ and $z$ and applying Lemma \ref{lem:p}, we obtain
\[
\int_{\R^{2d}} \varphi_{R,t}(s,y) \varphi_{R,t}(r,z) |y-z| ^{-\beta} dydz\le \frac C {\sigma_{R,t}^2} \int_{Q_R^2}  |x_1- x_2|^{-\beta} dx_1 dx_2.
\]
 The desired result follows from Lemma \ref{LemmaA.2} and \eqref{variance}.
\end{proof}

\begin{lemma}\label{lem:m} Let $m(t_{\alpha},R) $ be as in \eqref{m}. Then, there is a constant $C_{t,\beta,d}$ depending on $t$, $\beta$ and $t$, such that for $R\geq 1$ and $\e\leq 1$  
 \begin{align*}
    m(t_{\alpha},R) \geq C_{t,\beta,d}\e^{\alpha}.
\end{align*}
\end{lemma}
\begin{proof}
Integrating in $y$ using semigroup property and then applying the Fourier transform we get
\begin{align*}
m(t_{\alpha}, R) &=\frac{1}{
    \sigma_{R,t}^2}\int_0
    ^{\e^{\alpha}} \int_{\bm{Q}_R^2}  \int_{\R^{2d}} \bm{p}_{s}(x_1-y)\bm{p}_{s}(y-\ytilde-x_2)|\ytilde|^{-\beta} dy  d\ytilde dx_1dx_2 ds  \\
   & =\frac{1}{
    \sigma_{R,t}^2}\int_0^{\e^{\alpha}}\int_{\bm{Q}_R^2} \int_{\R^{d}} \bm{p}_{2s}(\ytilde-x_1+x_2)|\ytilde|^{-\beta} d\ytilde dx_1dx_2 ds \\
&=\frac{1}{
    \sigma_{R,t}^2}\int_0^{\e^{\alpha}}\int_{\bm{Q}_R^2}\int_{\R^{d}} e^{- i (x_1-x_2)\cdot \xi -s|\xi|^2} |\xi|^{\beta-d}   d\xi dx_1 dx_2 ds\\
    &=\frac{1}{
    \sigma_{R,t}^2}\int_0^{\e^{\alpha}}\int_{\R^{d}} e^{-s|\xi|^2} \left| \int_{\bm{Q}_R}e^{- i x\cdot \xi }dx\right|^2 |
    \xi|^{\beta-d}d\xi ds.
  \end{align*}  
   Integrating the variable $x$ and then the variable $s$ yields
   \begin{align*}
   m(t_{\alpha}, R) &=\frac{1}{
    \sigma_{R,t}^2}\int_{\R^{d}} \int_0^{\e^{\alpha}} e^{-s|\xi|^2} \left( \prod_{j=1}^d \frac{\sin(\xi_j R)}{\xi_j}\right)^2 |\xi|^{\beta-d}ds d\xi  \\
&= \frac{R^{2d}}{
    \sigma_{R,t}^2}\int_{\R^d} \frac{1-e^{-\epsilon^{\alpha}|\xi|^2}}{|\xi|^2} \left( \prod_{j=1}^d \frac{\sin(\xi_j R)}{\xi_jR}\right)^2 |\xi|^{\beta-d}d\xi .
  \end{align*}
 After the change of variable $\eta_j=R\xi_j$, we can rewrite the above integral as \begin{align*}
      m(t_\alpha,R)=\frac{R^{2d-\beta}\e^{\alpha}}{\sigma_{R,t}^2}  \int_{\R^d} \frac{1-e^{-\epsilon^{\alpha}|\eta|^2/R^2}}{\epsilon^{\alpha}|\eta|^2/R^2}\left( \prod_{j=1}^d \frac{\sin(\eta_j)}{\eta_j}\right)^2 |\eta|^{\beta-d}d\eta.
  \end{align*}
  Using the estimate $(1-e^{-x})/x \geq 1-x$  and $\left|\sin x/x\right| \geq c>0$ when $0\leq x<1$ together with \eqref{variance}, we get for $R\geq 1$ and $\e\leq 1$, 
\begin{align*}
 m(t_{\alpha},R) & \geq \frac{R^{2d-\beta}\e^{\alpha}}{\sigma_{R,t}^2}  \int_{\bm{B}(0,1)} \left({1-\epsilon^{\alpha}\frac{|\eta|^2}{R^2}} \right)\left( \prod_{j=1}^d \frac{\sin(\eta_j)}{\eta_j}\right)^2  |\eta|^{\beta-d}d\eta \\
 &\geq  c^{2d}\frac{R^{2d-\beta}\e^{\alpha}}{\sigma_{R,t}^2}  \int_{\bm{B}(0,1)} \left({1-\epsilon^{\alpha}\frac{|\eta|^2}{R^2}} \right)  |\eta|^{\beta-d}d\eta \\
    &\geq C_{t,\beta,d}\e^{\alpha}\int_0^1 \left({1-\epsilon^{\alpha}\frac{r^2}{R^2}} \right) r^{\beta-1}dr\\
    &=C_{t,\beta,d}e^{\alpha} \left(\frac{1}{\beta}-\frac{\e^{\alpha}}{(\beta+2)R^2}\right) \geq C_{t,\beta,d}e^{\alpha}.
\end{align*}
This completes the proof.
\end{proof}

\begin{lemma}\label{lem:E}
Let $ E_{R,t}(s_1,s_2,\tau) $ be as defined in \eqref{ert}. Then, there is a constant $C_{t,\beta,d}$ depending on $t$, $\beta$ and $d$, such that 
 \begin{align*}
    E_{R,t}(s_1,s_2,\tau) \leq C_{t,\beta,d} R^{4d-3\beta}.
\end{align*}
\end{lemma}

\begin{proof}
After applying the change of variables $\theta=\xi, \tilde{\theta}=\xi-\xitilde, \eta_1=y_1, \eta_2=y_2, \tilde{\eta}_1=y_1-\ytilde_1, \tilde{\eta}_2=y_2-\ytilde_2$, we can rewrite $E_{R,t}$ as follows
\begin{align*} 
 E_{R,t}(s_1,s_2,\tau) =\int_{\bm{Q}_R^4}\int_{\R^{6d}} & \bm{p}_{t-\tau}(x_1-\theta)\bm{p}_{t-\tau}(\theta-\tilde{\theta}-x_2)\bm{p}_{t-s_1}(\eta_1-\tilde{\eta}_1-x_3) \\
 & \times \bm{p}_{t-s_2}(\eta_2-\tilde{\eta}_2-x_4)\bm{p}_{\tau-s_1}(\theta-\eta_1)   \bm{p}_{\tau-s_2}(\theta-\tilde{\theta}-\eta_2)  \\
 & \times 
  |\tilde{\theta}|^{-\beta}\nonumber |\tilde{\eta}_1|^{-\beta}|\tilde{\eta}_2|^{-\beta} d\theta d\tilde{\theta}d\eta_1 d\eta_2 d\tilde{\eta_1} d\tilde{\eta_2} dx_1dx_2dx_3dx_4. \nonumber
\end{align*}
Integrating in $\eta_1$ and $\eta_2$ using semigroup property, and using Lemma~\ref{lem:p} for the integrals in $\tilde{\eta}_1$, $\tilde{\eta}_2$, we obtain \begin{align*}
  E_{R,t}(s_1,s_2,\tau) \leq  C_{d,\beta}\int_{\bm{Q}_R^4}\int_{\R^{2d}}& \bm{p}_{t-\tau}(x_1-\theta)\bm{p}_{t-\tau}(\theta-\tilde{\theta}-x_2) |x_3-\theta|^{-\beta} |x_4-\theta+\tilde{\theta}|^{-\beta} \\
  &  \times   |\tilde{\theta}|^{-\beta}d\theta d\tilde{\theta} dx_1dx_2dx_3dx_4  .
\end{align*}
Using Lemma \ref{LemmaA.2} we can estimate the integrals in $x_3$ and $x_4$ to get 
\[
  E_{R,t}(s_1,s_2,\tau) \leq R^{2d-2\beta}C_{d,\beta}\int_{\bm{Q}_R^2}\int_{\R^{2d}} \bm{p}_{t-\tau}(x_1-\theta)\bm{p}_{t-\tau}(\theta-\tilde{\theta}-x_2)  |\tilde{\theta}|^{-\beta} d\theta d\tilde{\theta} dx_1dx_2 .
\]
Finally, applying  again Lemma  \ref{lem:p}  twice we get 
\[
  E_{R,t}(s_1,s_2,\tau) \leq C_{d,\beta} R^{2d-2\beta} \int_{Q_R^2}  |x_1-x_2| ^{-\beta} dx_1 dx_2,
\]
which allows us to complete  the proof.
\end{proof}

\begin{lemma}\label{lem:conditionnegativemoments} 
       Let $F$ be a nonnegative random variable. Then $\ex{F^{-p}}<\infty$ for all $p \geq 2$ if and only if for all $q \geq 2$, there exists $C=C(q)>0$ and $\e_0=\e_0(q)>0$ such that $\pr{F<\e}\leq C \e^q$ for all $\e\leq \e_0$. 
\end{lemma}

{}

\end{document}